\documentclass[11pt]{article}
\usepackage{mathrsfs}
\usepackage{amsfonts}
\usepackage{pifont}
\usepackage{bbding}
\usepackage{latexsym}
\usepackage{mathrsfs}
\usepackage{amsmath}
\usepackage{amsthm}
\usepackage{amssymb}
\usepackage{graphicx}

\usepackage{amssymb, latexsym}
\usepackage{color,amsmath,amssymb, amsfonts,amstext,amsthm}

\usepackage{epsfig, graphicx, graphics}
\usepackage{longtable}
\usepackage{subfigure}

\newcommand{\intd}{\textrm{d}}

\newcommand{\dps}{\displaystyle}
\usepackage{cite}

\newtheorem{theorem}{\indent Theorem}[section]

\newtheorem{definition}{\indent Definition}[section]
\newtheorem{lemma}{\indent Lemma}[section]
\newtheorem{remark}{\indent Remark}[section]
\numberwithin{equation}{section}

\title{Approximation of the inertial manifold\\ for a nonlocal dynamical system}

\author{Xingjie Yan\\ Department of Mathematics, China University of Mining and Technology,\\ Xuzhou, Jiangsu  221008,  China\\
Jinchun He\\
School of Mathematics and Statistics,\\
Huazhong University of Science and Technology, Wuhan,   430074, China\\
Jinqiao Duan\\Department of Applied Mathematics, Illinois Institute of Technology,\\ Chicago, Illinois, 60616, USA
}

\begin{document}

\date{January 30, 2014}

\maketitle

\footnote{This work was partly supported by the Fundamental Research Funds for the Central Universities Grants (2010QNA40, 2012LWN56, 2010QNB20), the Tianyuan Foundation of China grant (11126306), and the National Nature Science Foundation of Jiangsu Province grant (BK20130170) and the National Nature Science Foundation of China grant (11271364).  A part of this work was done while the authors were visiting the Institute for Pure and Applied Mathematics (IPAM), UCLA, Los Angeles, USA.\\
Corresponding author: J. He  (taoismnature@mail.hust.edu.cn)  }

\begin{abstract}
We consider  inertial manifolds and their approximation for a class of partial differential equations with  a nonlocal Laplacian operator $-(-\Delta)^{\frac{\alpha}{2}}$, with $0<\alpha<2$.  The nonlocal or fractional Laplacian operator represents  an anomalous diffusion  effect.   We first establish the existence of an inertial manifold and highlight the influence of the parameter  $\alpha$. Then we approximate the inertial manifold when a small normal diffusion   $\varepsilon \Delta$ (with $\varepsilon \in (0, 1)$) enters the system, and obtain the estimate for the Hausdorff semi-distance between the inertial manifolds with and without normal diffusion.

\medskip

\textbf{Keywords}:  Anomalous and normal diffusion; Inertial manifolds; Nonlocal dynamical systems; Partial differential equations with a fractional Laplacian operator

\end{abstract}

\section{Introduction}

Nonlocal operators appear in complex systems,  such as anomalous diffusion and geophysical flows \cite{18, MS2012, MK2004},  a thin obstacle problem \cite{14}, finance \cite{29}, and stratified materials \cite{30}.
A special but important   nonlocal operator  is the fractional Laplacian operator arising in non-Gaussian stochastic systems.
For a  stochastic differential system with  a symmetric $\alpha$-stable L\'evy motion (a non-Gaussian stochastic process) $L_t^\alpha$ for $\alpha \in (0, 2)$,
$$ dX_t =b(X_t)dt + dL^\alpha_t, \;\; X_0 =x, $$
the corresponding Fokker-Planck equation contains the fractional Laplacian operator $(-\Delta)^{\frac{\alpha}{2}}$. When the drift term  $b$ in the above  stochastic differential system depends on the probability distribution of the system state, the Fokker-Planck equation becomes a nonlinear, nonlocal partial differential equation \cite{Applebaum}.

For nonlocal  partial differential equations with the fractional Laplacian operator $(-\Delta)^{\frac{\alpha}{2}}$, $\alpha \in (0, 2)$, there are     recent works about  modeling techniques, well-posedness and regularity of solutions; see, for example \cite{13, 14,15,16,18,26,28,CMN2012, 12}. It is desirable to further study dynamical behaviors of such nonlocal systems.



In the present paper, we consider the inertial manifold and its approximation for a  system described by a nonlocal partial differential equation.
An inertial manifold is a Lipschitz manifold that captures asymptotic long time dynamics of the   system evolution \cite{1, 2, 7}.
The Lyapunov-Perron method is often used to study an inertial manifold. With this method, the system is converted  into an integral equation and the inertial manifold is constructed as the graph of a unique fixed point of a corresponding mapping. The method is also used to construct the inertial manifolds of stochastic partial differential equations  \cite{20}.
  A spectral gap condition is sufficient to guarantee the existence of the fixed point. This spectral gap condition may be understood as a relationship  of  the gap between eigenvalues of the linear operator in the system, with the Lipschitz constant of the nonlinearity. See the inequality (\ref{gap}) in the next section.

We consider a nonlocal dynamical system containing anomalous diffusion $(-\Delta)^{\frac{\alpha}{2}}$, with or without small normal diffusion $ \varepsilon \Delta$ for $\varepsilon \in (0, 1)$. We prove the existence of an inertial manifold and then consider its approximation when $\varepsilon$ is sufficiently small.


The paper is arranged as follows. In   section 2, we review a basic theory of inertial manifolds. In   section 3, we   prove the existence of the inertial manifold  of a nonlocal   system  with both anomalous and normal diffusion. We point out that the interaction between the normal diffusion and the anomalous diffusion in the nonlocal operator plays a significant role in the existence of the inertial manifold. Section 4 is  devoted to the existence of the inertial manifold  when the normal diffusion   is absent; however, in this case we do not  have the existence of the inertial manifold  when the parameter $\alpha$ in the anomalous diffusion is less than $1$, since the spectral gap condition does not hold. Finally, in section 5, we derive an asymptotic approximation of the inertial manifold when the normal diffusion is small enough.


\section{Preliminaries}
\label{prelim}

In this section, we recall   basic facts about inertial manifolds of infinite dimensional dynamical systems  \cite{1,2,3}. We consider an evolution equation in a Hilbert space $H$ with norm $|\cdot|$ and scalar product $\langle\cdot,\cdot\rangle$
\begin{equation}\label{e2.1}
\frac{du}{dt}+Au = f(u),~u\in H,~u|_{t=0}=u_0,
\end{equation}
where $A$ is a closed linear operator on $H$ and $f$ is a nonlinear mapping. We make the following assumptions.

{\bf Assumption I}.  The linear operator $A$ is a positive definite, self-adjoint operator with   discrete spectrum. More specifically, the eigenvalues  may be arranged as follows
\begin{equation}0 < \lambda_1 < \lambda_2 \leq \cdot\cdot\cdot\leq \lambda_j
\leq \cdot\cdot\cdot,~~~\lambda_j \rightarrow + \infty
~\mathrm{as}~j \rightarrow \infty.
\end{equation}
For example, if a linear closed operator $A$   is positive and self-adjoint   with the compact inverse, then it has   discrete spectrum.

Let $V$ be a subspace of $H$ with norm $\|\cdot\|$ and scalar product $\ll\cdot,\cdot\gg$.

{\bf Assumption II}.  The nonlinear mapping $f: V\longrightarrow H$ is globally Lipschitz continuous
\begin{equation}
|f(u)-f(v)|\leq l_f\|u-v\|,~\forall~ u,v\in V,
\end{equation}
where $l_f$ is the Lipschitz constant.

It is well-known that if the operator $-A$ is the infinitesimal generator of a $C_0-$semigroup, the problem (\ref{e2.1}) is well-posed, that is for every $u_0\in H$, there exists a unique global solution $u(t)$,   and, in fact, $u\in C(0,T;H)\cap L^2(0,T;H^1_0(\Omega))$ for all $T>0$; see \cite{4}. The solution of the system (\ref{e2.1}) can be defined by the formula of variation of constant
\begin{align*}
u(t)=e^{-At}u_0 + \int^t_0e^{-A(t-s)}f(u(s))ds.
\end{align*}
If we define $T(t):u_0\longrightarrow u(t)$, then $T(t)$ is continuous and satisfies semigroup properties:
\begin{enumerate}
\item[(i)]$T(0)=I  \textrm{(identity~operator)}$;
\item[(ii)]$T(t)T(s)=T(t+s)$, $t,s\geq0$.
\end{enumerate}

Let $P_n$ a projection operator with rank $n$ from $H$ to $P_nH$, and then we define $Q_n=I-P_n$.

\begin{definition}

We say that a manifold $\mathcal {M}$ in space $H$ is an inertial manifold for the  system (\ref{e2.1}) or the corresponding semigroup $T(t)$ if
\begin{enumerate}
\item[(i)] the manifold $\mathcal {M}$ is invariant under $T(t)$, i.e., $T(t)\mathcal {M}=\mathcal {M}$, $t\geq0$;
\item[(ii)] it can be represented as the graph of a Lipschitz continuous function $\Phi:P_nH\longrightarrow Q_nH$, i.e.,
$$\mathcal {M}=\{u+\Phi(u),u\in P_nH\};$$
\item[(iii)] it  possesses exponential tracking property, i.e., there exist positive constants $\eta$ and $\beta$ such that, for every $u\in H$, there is a $v\in\mathcal {M}$ such that
    $$|T(t)u-T(t)v|\leq \eta e^{-\beta t}|u-v|.$$
\end{enumerate}
\end{definition}

Under certain further assumptions, for a large enough $N$, such that $\lambda_{N+1}-\lambda_N>2l_f$, an inertial manifold for system (\ref{e2.1}) can be realized as the graph of a function $\Phi:P_N H\longrightarrow Q_N H$. The inertial form
\begin{equation}
\frac{dp}{dt}+ Ap=P_N f(p+\Phi(p)),~p\in P_N H,
\end{equation}
captures  all long-time behaviors of (\ref{e2.1}). The dimension of $\mathcal {M}$ is the dimension of $P_N H$, i.e., $$dim(\mathcal {M})=dim(P_NH).$$

We assume that there exists an exponential dichotomy (see \cite{6}):

\begin{enumerate}
\item[(i)]$\|e^{-tA}P_n\|_{\mathcal {L}(H,H)}\leq K_1e^{-\lambda_nt},~\forall t\leq0$,
\item[(ii)]$\|e^{-tA}P_n\|_{\mathcal {L}(V,H)}\leq K_1\lambda_n^se^{-\lambda_nt},~\forall t\leq0$,
\item[(iii)]$\|e^{-tA}(I-P_n)\|_{\mathcal {L}(H,H)}\leq K_2e^{-\lambda_{n+1}t},~\forall t\geq0$,
\item[(iv)]$\|e^{-tA}(I-P_n)\|_{\mathcal {L}(V,H)}\leq K_2(t^{-s}+\lambda_{n+1}^s)e^{-\lambda_{n+1}t},~\forall t\geq0$,
\end{enumerate}
where  constants $K_1,~K_2\geq1,$ $0\leq s<1$, and $\|\cdot\|_{\mathcal {L}(V,H)}$ is   operator norm.

If a spectral gap condition holds, i.e.,   there is an $N\in\mathbb{N}$  such that
\begin{equation}
\lambda_{N+1}-\lambda_N>2l_f,  \label{gap}
\end{equation}
then, we can choose a $\sigma$ such that
\begin{equation}
\lambda_N + 2l_fK_1<\sigma<\lambda_{N+1}-2l_fK_1K_2.
\end{equation}
This $\sigma$ is used to define the Banach space
$$\mathcal {F}_\sigma=\{\varphi\in C((-\infty,0],H)~|~\|\varphi\|_{\mathcal {F}_\sigma}=\sup_{t\leq0}e^{\sigma t}|\varphi(t)|<\infty\}.$$
A trajectory on the inertial manifold can be found as the fixed point $\varphi=\varphi(p)$ of a mapping $J(\cdot,p):\mathcal {F}_\sigma\longrightarrow\mathcal {F}_\sigma$ defined  by
\begin{equation}\label{e2.6}
J(\varphi,p)(t)=e^{-tA}p-\int^0_te^{-(t-s)A}P_Nf(\varphi(s))ds+\int^t_{-\infty}e^{-(t-s)A}(I-P_N)f(\varphi(s))ds.
\end{equation}
The inertial manifold $\mathcal {M}$ is the graph of  $\Phi:P_N H\longrightarrow(I-P_N)H$, which is defined in terms of  the fixed point $\varphi$  of  (\ref{e2.6}) as follows
\begin{align}\label{e2.7}
\Phi(p)=(I-P_N)\varphi(p)(0)&=(I-P_N)(\int^0_{-\infty}e^{sA}(I-P_N)f(\varphi(s))ds)\nonumber\\&=\int^0_{-\infty}e^{sA}(I-P_N)f(\varphi(s))ds,~\forall~p\in P_N H.
\end{align}
Note that $\varphi(p)(0)=p+\Phi(p)$.

The spectral gap condition is used to ensure not only that $J$ has a fixed point by contraction mapping principle, but also that the resulting manifold is exponentially tracking, i.e., there exist   positive constants $\eta$ and $\beta$ such that, for   $u\in H$, there is $v_0\in\mathcal {M}$ such that
    $$|T(t)u-T(t)v_0|\leq \eta e^{-\beta t}|u-v_0|.$$

The conclusion mentioned above is summarized in the following theorem (see \cite{7,1})

\begin{theorem}\label{t2.1}
Let the Assumptions I and II hold. If there is an $N\in\mathbb{N}$, such that the following spectral gap condition is satisfied
\begin{equation}
\lambda_{N+1}-\lambda_N>2l_f,
\end{equation}
then there exists an $N$-dimensional inertial manifold $\mathcal {M}$, which is the graph of a Lipschitz continuous function $\Phi$ satisfying (\ref{e2.7}).
\end{theorem}

\begin{remark}
If $f$ is $C^{k}$, then $\mathcal {M}$ is $C^{k}$.

\end{remark}

\begin{remark}
If the system (\ref{e2.1}) is dissipative, i.e., it possesses an absorbing set, then   $f$ only needs to be locally Lipschitz continuous.
\end{remark}

\section{Inertial manifold for a nonlocal system with both anomalous and normal diffusion}
\label{usual}

Now we consider the existence of the inertial manifold for a  nonlocal evolution equation with both fractional and usual Laplacian
(i.e.,   both   anomalous and normal diffusion).  We recall the definition of the fractional Laplacian operator.
\begin{definition}
For $u\in C^\infty_0(\mathbb{R}^n)$ and  $\alpha\in(0,2)$, define
$$(-\Delta)^{\frac{\alpha}{2}} u=C(n,\alpha)P.V.\int_{\mathbb{R}^n}\frac{u(x)-u(y)}{|x-y|^{n+2\alpha}}dy,$$
where the principle value (P.V.) is taken as the limit of the integral over $\mathbb{R}^n\backslash B_\epsilon(x)$ as $\epsilon\longrightarrow0$, with $B_\epsilon(x)$   the ball  of radius $\epsilon$ centered at $x$, and $$C(n,\alpha)=\frac{\alpha}{2^{1-\alpha}\pi^{\frac{n}{2}}}\frac{\Gamma(\frac{n+\alpha}{2})}{\Gamma(1-\frac{\alpha}{2})}.$$
Here $\Gamma$ is the Gamma function defined by $\Gamma(\lambda)=\int^\infty_0t^{\lambda-1}e^{-t}dt$ for every $\lambda>0$; for more information see\cite{25,26}. In this paper, $n=1$ and the usual local Laplacian operator is $\Delta = \partial_{xx}$.
\end{definition}

We consider the existence of the inertial manifold  $\mathcal{M}_\varepsilon$ for the following system
\begin{align}\label{3.1}
\begin{cases}
\dfrac{du}{dt} - \varepsilon\Delta u +(-\Delta)^{\frac{\alpha}{2}} u+f(u)=g(x)
,&in\  \Omega=(-\pi,\pi).\\
u\big|_{\Omega^c}=0.\\
u(x,0)=u_0,  &\  x\in \Omega,
\end{cases}
\end{align}
where $0<\varepsilon<<1$, and $\Omega^c=\mathbb{R}\backslash\Omega$. Let $H=L^2(\Omega)$.  The nonlocal Laplacian operator $(-\Delta)^{\frac{\alpha}{2}}$ is defined on $H$.   Assume that the nonlinear function $f:H\longrightarrow H$ is locally Lipschitz, i.e., for every $B\subset H$
\begin{equation}\label{eee3.2}
|f(u)-f(v)|\leq l_f|u-v|,~\forall~ u,v\in B.
\end{equation}

{\bf Assumption III}.  The nonlinear function $f$ satisfies the following condition
$$C|s|^{p}-C\leq f(s)s\leq C|s|^{p}+C,\;\;  s \in \mathbb{R},  $$
 for some $p\geq2$.

We recall the eigenvalues of the nonlocal operator $(-\Delta)^{\frac{\alpha}{2}}$ in $H=L^2(\Omega)$.

\begin{lemma}(\cite{10})\label{l3.1}
The eigenvalues of the following spectral problem
\begin{equation}
(-\Delta)^{\frac{\alpha}{2}}\varphi(x)=\lambda\varphi(x),~x\in\Omega,
\end{equation}
where $\varphi(x)\in L^2(\Omega)$ is extended to $\mathbb{R}$ by $0$, are
\begin{equation}
\lambda_n=\Big(\frac{n}{2}-\frac{(2-\alpha)}{8}\Big)^\alpha+o(\frac{1}{n}).
\end{equation}
Moreover,
$$0 < \lambda_1 < \lambda_2 \leq \cdot\cdot\cdot\leq \lambda_n
\leq \cdot\cdot\cdot, \; \mbox{for}~ n=1, 2, \cdots.
$$
Moreover, the corresponding eigenfunctions $\varphi_n$ form a complete orthonormal basis in $L^2(\Omega)$.
\end{lemma}

Then we have the following result on the well-posedness    for the system (\ref{3.1}).   For some related results, see \cite{5}.
Denote by     $C$ a general positive constant which may be different in different places.

\begin{theorem}\label{t3.1} (Well-posedness)
Assume that $g$ is in $L^2(\Omega)$, $f$ satisfies the condition (\ref{eee3.2}), and that Assumption III hold. Then there exists a unique solution $u(t)\in C(0,T;H)$ (for every $T>0$)  for  the system (\ref{3.1}). The solution is given by the formula of   variation of constants
$$u(t)=e^{-A_\varepsilon t}u_0 - \int^t_0e^{-A_\varepsilon(t-s)}(f(u(s))-g(x))ds,$$
where $A_\varepsilon=- \varepsilon\Delta +(-\Delta)^{\frac{\alpha}{2}} $.
\end{theorem}

\begin{proof}
Recall that  the eigenvalues of the local Laplacian operator $- \Delta$ in $L^2(\Omega)$ with domain $H^2(\Omega)\cap H^1_0(\Omega)$ are $ n^2$, $n\in\mathbb{N}$, and the corresponding eigenfunctions are  $\omega_n(x)=\sin(nx)$ which form a complete orthonormal basis of $L^2(\Omega)$. By Lemma \ref{l3.1} and \cite{31}, we know that the eigenvalues of $A_\varepsilon$  are $\lambda_n=\varepsilon n^2+(\frac{n}{2}-\frac{(2-\alpha)}{8})^\alpha+o(\frac{1}{n})$ $n\in\mathbb{Z}$, $n\geq1$, and they satisfy $0 < \lambda_1 < \lambda_2 \leq \cdot\cdot\cdot\leq \lambda_n
\leq \cdot\cdot\cdot$. The corresponding eigenfunctions form a complete orthonormal basis in $L^2(\Omega)$. Hence, we   conclude that the operator $-A_\varepsilon=\varepsilon\Delta  -(-\Delta)^{\frac{\alpha}{2}} $ is dissipative in $L^2(\Omega)$, that is, $<-A_\varepsilon u,u>\leq0$. Thus, $-A_\varepsilon$ is a infinitesimal generator of an analytic semigroup \cite{8}. As $f$ is local Lipschitz continuous, we obtain the existence and uniqueness of solution $u(t)\in C(0,T;H)\cap L^2(0,\tau(u_0);H^1_0(\Omega))$ for   some $\tau=\tau(u_0)>0$; see \cite{4}. In fact, the solution is given by the formula of   variation of constants. Now,  we prove that the solution   exists globally.

Multiplying $u$ on   both sides of the equation (\ref{3.1}) and integrating over $\Omega$, we have
\begin{eqnarray}
\frac12\frac{\intd}{\intd t}|u|^2+\varepsilon\|u\|^2 +\int_\Omega f(u)uds&=-&\int_\Omega u(-\Delta)^{\frac{\alpha}{2}}u\intd x+\int_\Omega g(x)u \intd x.
\end{eqnarray}
Furthermore,
\begin{eqnarray}
\frac12\frac{\intd}{\intd t}|u|^2+\varepsilon\|u\|^2 +C\int_\Omega|u|^pdx\le -\int_\Omega\int_\Omega (\mathcal{D}^{\ast}(u))^2\intd y \intd x+C\int_\Omega u^2 \intd x+C|g(x)|^2 +C|\Omega|,
\end{eqnarray}
where $|\Omega|$ denote the measure of $\Omega$.
Using the nonlocal Poincar\'e $|u|^2\le C\int_\Omega\int_\Omega (\mathcal{D}^{\ast}(u))^2\intd y \intd x$ (see \cite{DGLZ2011}), and local Poincar\'e inequality, we get
\begin{equation}
\frac{\intd}{\intd t}|u|^2+C|u|^2 \le C|g(x)|^2+C|\Omega|.
\end{equation}
By uniform Gronwall's inequality, we have
\begin{equation}
|u|^2 \le |u_0|^2\textrm{e}^{-Ct}+C|g(x)|^2+C|\Omega|.
\end{equation}
which implies that $\displaystyle{\sup_{0\le t<+\infty}|u(x,\,t)|<+\infty}$. Thus the solution exists for all time. The proof is complete.
\end{proof}

By Theorem \ref{t3.1}, we can define a semigroup by $T(t):u_0\longrightarrow u(t)$, where $u(t)$ is the solution of   (\ref{3.1}) and the operator $-A_\varepsilon$ is the infinitesimal generator of $T(t)$. Furthermore, for the semigroup $\{T(t)\}_{t\geq0}$, there exists a bounded absorbing set $B_\varepsilon$ in $H$, i.e., for every bounded set $B\subset H$, we can find a constant $t_0=t_0(B)>0$, such that when $t\geq t_0$, $T(t)B\subset B_\varepsilon$.

We are ready to present the   result on the existence of the inertial manifold $\mathcal {M}_\varepsilon$.

\begin{theorem}\label{t3.2} (Inertial manifold $\mathcal {M}_\varepsilon$) \\
Assume that $g$ belongs to $L^2(\Omega)$, the nonlinearity $f$  is local Lipschitz continuous, and   Assumption III holds. Then there exists an $N-$dimensional inertial manifold $\mathcal {M}_\varepsilon$ for the system (\ref{3.1}), as a graph of a Lipschitz continuous function $\Phi$ from $P_NH$ to $(I-P_N)H$.
\end{theorem}

\begin{proof}
Both $-\Delta$ and $(-\Delta)^{\frac{\alpha}{2}}$, $\alpha\in(0,2)$ are positive, seif-adjoint operators in $L^2(\Omega)$. By Lemma \ref{3.1} and \cite{31},   the eigenvalues of the operator $A_\varepsilon=- \varepsilon\Delta  +(-\Delta)^{\frac{\alpha}{2}}$     are
\begin{equation}
\lambda_n=\varepsilon n^2+\dps{\left(\frac{n}{2}-\frac{(2-\alpha)}{8}\right)}^\alpha+o(\frac{1}{n}),
\end{equation}
which satisfy
$$0 < \lambda_1 < \lambda_2 \leq \cdot\cdot\cdot\leq \lambda_n
\leq \cdot\cdot\cdot,$$
and the corresponding eigenfunctions $\varphi_n$ form a complete orthonormal basis of $L^2(\Omega)$. The assumptions I and II are all satisfied. Hence, by Theorem \ref{t2.1}, we only need to verify the spectral gap condition. We do this in the  following cases:

{\bf Case 1, $\alpha=1$:} $\lambda_{n+1}-\lambda_{n}=\varepsilon(2n+1)+\frac{1}{2}+o(\frac{1}{n}).$ So the spectral gap condition is satisfied for fixed $\varepsilon$, i.e., we can find some $N$ such that $\lambda_{N+1}-\lambda_{N}\geq 2l_f$, while $n\geq N$, we obtain $\lambda_{n+1}-\lambda_{n}\geq\lambda_{N+1}-\lambda_{N}\geq 2l_f$.

Next, suppose that $G(n)=\lambda_{n+1}-\lambda_{n}=\varepsilon(2n+1)+ (\frac{n+1}{2}-\frac{2-\alpha}{8})^\alpha-(\frac{n}{2}-\frac{2-\alpha}{8})^\alpha+o(\frac{1}{n})$. Then
\begin{equation}\label{e3.6}
G^{'}(n) = 2\varepsilon+\frac{\alpha}{2}\Big[\big(\frac{n+1}{2}-\frac{2-\alpha}{8}\big)^{\alpha-1}-\big(\frac{n}{2}-\frac{2-\alpha}{8}\big)^{\alpha-1}\Big]+\gamma o(\frac{1}{n^2}).
\end{equation}
If $G^{'}(n)>0$, then the spectral gap becomes larger and larger, hence for $2l_f$, we can find suitable $N$ such that spectral gap condition hold.

{\bf Case 2, $1<\alpha<2$:} In this case, the term $\frac{\alpha}{2}[(\frac{n+1}{2}-\frac{2-\alpha}{8})^{\alpha-1}-(\frac{n}{2}-\frac{2-\alpha}{8})^{\alpha-1}]>0$, and $\gamma o(\frac{1}{n^2})$ becomes very small when $n$ large enough whenever the sign of $\gamma$. Hence we can choose big $n$, such that $G^{'}(n)>0$, so the gap becomes larger and larger. Hence for $2l_f$, we can choose big $N$, such that $n\geq N$, we obtain $\lambda_{n+1}-\lambda_{n}\geq\lambda_{N+1}-\lambda_{N}\geq 2l_f$.

{\bf Case 3, $0<\alpha<1$:} From (\ref{e3.6}), we know that in this situation $\frac{\alpha}{2}[(\frac{n+1}{2}-\frac{2-\alpha}{8})^{\alpha-1}-(\frac{n}{2}-\frac{2-\alpha}{8})^{\alpha-1}]<0.$ Our aim is to obtain $\lambda_{n+1}-\lambda_{n}$ larger and larger, so necessary, we need $G^{'}(n)>0$. By above analysis, if $\gamma o(\frac{1}{n^2})>0$, then we choose $-\frac{\alpha}{4}[(\frac{n+1}{2}-\frac{2-\alpha}{8})^{\alpha-1}-(\frac{n}{2}-\frac{2-\alpha}{8})^{\alpha-1}]-\frac{\gamma }{2} o(\frac{1}{n^2})<\varepsilon<<1$, such that $G^{'}(n)>0$, i.e., there exist $N$, such that $n\geq N$, we obtain $\lambda_{n+1}-\lambda_{n}\geq\lambda_{N+1}-\lambda_{N}\geq 2l_f$. If $\gamma o(\frac{1}{n^2})<0$, then we choose $-\frac{\alpha}{4}[(\frac{n+1}{2}-\frac{2-\alpha}{8})^{\alpha-1}-(\frac{n}{2}-\frac{2-\alpha}{8})^{\alpha-1}]+\frac{\gamma }{2} o(\frac{1}{n^2})<\varepsilon<<1$, such that $G^{'}(n)>0$, i.e., there exists $N$, such that $n\geq N$, we obtain $\lambda_{n+1}-\lambda_{n}\geq\lambda_{N+1}-\lambda_{N}\geq 2l_f$.

In conclusion, we see that in cases 1 and 2, for arbitrary $\varepsilon>0$ , there exists inertial manifold, but in case 3, the existence of inertial manifold under the choice of $\varepsilon$. The proof is complete.
\end{proof}

\section{Inertial manifold for a nonlocal   system   with only  anomalous diffusion }
\label{frac}

In this section, we consider the existence of the inertial manifold   $\mathcal {M}_0$ of following   equation with only anomalous diffusion
\begin{align}\label{4.1}
\begin{cases}
\dfrac{du}{dt} +(-\Delta)^{\frac{\alpha}{2}}+f(u)=g(x)
,&in\  \Omega=(-\pi,\pi).\\
u\big|_{\Omega^c}=0.\\
u(x,0)=u_0(x)\in H^{\frac{\alpha}{2}}_0(\Omega),&in\  x\in \Omega,
\end{cases}
\end{align}
where $\Omega^c=\mathbb{R}\backslash\Omega$.

 Assume that for any $B\subset H^{\frac{\alpha}{2}}_0(\Omega)$ function $f:B\longrightarrow H^{\frac{\alpha}{2}}_0(\Omega)$ is locally Lipschitz continuous
\begin{equation}\label{e3.2}
\|f(u)-f(v)\|_{H^{\frac{\alpha}{2}}_0}\leq l_f\|u-v\|_{H^{\frac{\alpha}{2}}_0},~\forall~ u,v\in B.
\end{equation}

We will prove the following Theorem.

\begin{theorem}\label{t4.1} (Inertial manifold $\mathcal {M}_0$)
Assume that $g$ belongs to $L^2(\Omega)$ and the nonlinearity of problem (\ref{4.1}) be local Lipschitz continuous (\ref{e3.2}). Also suppose that Assumption $III$ hold. Then there exists an $N-$dimensional inertial manifold $\mathcal {M}_0$ which is a graph of Lipschitz continuous function $\Phi$ from $P_NH^{\frac{\alpha}{2}}(\Omega)$ to $(I-P_N)H^{\frac{\alpha}{2}}(\Omega)$.
\end{theorem}

First, we begin with the existence and uniqueness of problem (\ref{4.1}), we will use semigroup method (\cite{8}). Denote $A_\alpha=-(-\Delta)^{\frac{\alpha}{2}}$ and $\|\ \|$ norm of space or operator.

\subsection{Some estimates on the nonlocal Laplacian}\label{semigroup-esti}

\begin{definition}
$A$ is a sectorial operator, if $A$ is dense defined for some $\phi\in(0,\,\frac{\pi}{2})$,
$M\ge 1$, and $a\in\mathbb{R}$, $S_{a,\,\phi}=\{\lambda|~\phi\le|\mbox{\textrm{arg}}(\lambda-a)|\le\pi,~\lambda\neq a\}
\subset\rho(A)$ and $\|(\lambda I-A)^{-1}\|\le M/|\lambda-a|$.
\end{definition}

\begin{lemma}\label{semiest}
The nonlocal Laplacian operator $A_{\alpha}$ is sectorial one, satisfying the estimates as follows
\begin{equation}
\|\textrm{e}^{tA_{\alpha}}\|_{L^2(\Omega)}\le C\textrm{e}^{-\delta t},~~~~\|A_{\alpha}\textrm{e}^{tA_{\alpha}}\|_{L^2(\Omega)}
\le \frac{C}{t}\textrm{e}^{-\delta t},
\end{equation}
where $C,\,\delta>0$ are constants independent of $t$.
\end{lemma}
\begin{proof}
By Lemma \ref{l3.1}, $A_{\alpha}$ is sectorial can be proved by definition. Set $\mu=\lambda t$, $\lambda>0$,
$$
\|\textrm{e}^{tA_{\alpha}}\|_{L^2(\Omega)}=\dps{\left\|\frac{1}{2\pi i}\int_{\Gamma}\textrm{e}^{\mu}(\frac{\mu}{t}-A_{\alpha})^{-1}
\frac{\intd\mu}{t}\right\|}_{L^2(\Omega)} \le \frac{M}{2\pi}\int_{\Gamma}\mid\textrm{e}^{\mu}\mid
\frac{\mid\intd\mu\mid}{\mid \mu\mid} \le C\textrm{e}^{-\delta t}.
$$
$$
\|A_{\alpha}\textrm{e}^{tA_{\alpha}}\|_{L^2(\Omega)}=\dps{\left\|\frac{1}{2\pi i}A_{\alpha}\int_{\Gamma}\textrm{e}^{\mu}(\frac{\mu}{t}-A_{\alpha})^{-1}
\frac{\intd\mu}{t}\right\|}_{L^2(\Omega)} \le \frac{1}{2\pi}\frac{M}{\delta}\int_{\Gamma}\mid\textrm{e}^{\mu}\mid
\frac{\mid\intd\mu\mid}{\mid \mu\mid}\frac{1}{t} \le \frac{C}{t}\textrm{e}^{-\delta t}.
$$
The proof is complete.
\end{proof}

\begin{definition}\label{indexoperator}
$A$ is a sectorial operator. if $\Re \sigma(-A)>0$, then for every $\beta>0$,
$A^{-\beta}=\dps{\frac{1}{\Gamma(\beta)}\int_0^{\infty}t^{\beta-1}\textrm{e}^{tA}\intd t}$.
Moreover, $A^{\beta}=(A^{-\beta})^{-1}$ and $A^0=I$.
\end{definition}

\begin{lemma}\label{semiregu}
$A$ is a sectorial operator, $\Re \sigma(-A)>\delta>0$. $\forall\,\beta\ge 0$, such that
$\exists\,C(\beta)<\infty$, $\forall\,t>0$,
\begin{equation}
\|A^{\beta}\textrm{e}^{tA}\|\le C(\beta) t^{-\beta}\textrm{e}^{-\delta t}.
\end{equation}
\end{lemma}
\begin{proof}
$\forall$ $m=1,\,2,\,\cdots$,
\begin{equation}
\|A^m\textrm{e}^{tA}\|=\left\|\left(A\textrm{e}^{\frac{tA}{m}}\right)^m\right\|
\le (Cm)^m t^{-m}\textrm{e}^{-\delta t}.
\end{equation}
For $0<\beta<1$, $t>0$
\begin{eqnarray*}
\|A^{\beta}\textrm{e}^{tA}\|&=&\|A^{\beta-1}A\textrm{e}^{tA}\|=\|A^{-(1-\beta)}A\textrm{e}^{tA}\| \\
&=&\dps{\left\|\frac{1}{\Gamma(1-\beta)}\int_0^{\infty}\tau^{1-\beta-1}A\textrm{e}^{-A(t+\tau)}
\intd \tau\right\|}\le \frac{1}{\Gamma(1-\beta)}\int_0^{\infty}\tau^{-\beta}
\|A\textrm{e}^{-A(t+\tau)}\| \intd \tau \\
&\le& \frac{C}{\Gamma(1-\beta)}\int_0^{\infty}\tau^{-\beta}(t+\tau)^{-1}\textrm{e}^{-\delta(t+\tau)}
\intd \tau=C\Gamma(\beta)t^{-\beta}\textrm{e}^{-\delta t}.
\end{eqnarray*}
Hence, $\forall$ $\beta\ge 0$, $\|A^{\beta}\textrm{e}^{tA}\|\le C(\beta) t^{-\beta}\textrm{e}^{-\delta t}$. The proof is complete.
\end{proof}

\begin{lemma}\label{embetoHs}
$\textrm{Dom}(A_{\alpha}) \hookrightarrow H_0^{\frac{\alpha}{2}}(\Omega)$.
\end{lemma}
\begin{proof}

In paper \cite{CMN2012}, we know $\textrm{Dom}(A_{\alpha}) \subset H_0^{\frac{\alpha}{2}}(\Omega)$. According to the embedding
results and the nonlocal calculus in paper \cite{DGLZ2011}, we have by the H\"{o}lder inequality
\begin{eqnarray*}
\|u\|^2_{H^{\frac{\alpha}{2}}(\Omega)}&=&\|u\|^2_{L^2(\Omega)}+|u|^2_{H^{\frac{\alpha}{2}}(\Omega)} \\
&\le& \|u\|^2_{L^2(\Omega)}+\frac12\int_\Omega\int_\Omega(\mathcal{D}^{\ast}(u)(x,\,y))^2 \intd x \intd y
+C\varepsilon^{-(1+\alpha)}\|u\|^2_{L^2(D)}\\
&=& C(\Omega,\,\varepsilon,\,\alpha)\|u\|^2_{L^2(\Omega)}-\frac12 \langle A_{\alpha}u,\,u \rangle \le C\|u\|^2_{L^2(\Omega)}
+C\|A_{\alpha}u\|_{L^2(\Omega)}\|u\|_{L^2(\Omega)} \\
&=& C\|u\|_{L^2(\Omega)}(\|u\|_{L^2(\Omega)}+\|A_{\alpha}u\|_{L^2(\Omega)}).
\end{eqnarray*}
Through the nonlocal Poincar\`{e} inequality\cite{DGLZ2011}, we get
\begin{eqnarray*}
\|u\|^2_{H^{\frac{\alpha}{2}}(\Omega)}&\le& C\|u\|_{H^{\frac{\alpha}{2}}(\Omega)}(\|u\|_{L^2(\Omega)}+\|A_{\alpha}u\|_{L^2(\Omega)}).
\end{eqnarray*}
Therefore, $\|u\|_{H^{\frac{\alpha}{2}}(\Omega)}\le C(\|u\|_{L^2(\Omega)}+\|A_{\alpha}u\|_{L^2(\Omega)})$. The proof is complete.
\end{proof}

\begin{lemma}\label{fracembe}
$\textrm{Dom}(A_{\alpha}^{\beta}) \hookrightarrow H^{\frac{\alpha}{2}}(\Omega)$, $\frac12<\beta<1$.
\end{lemma}
\begin{proof}
Through the Lemma $\ref{embetoHs}$, we have $\textrm{Dom}(A_{\alpha}^{\beta})
 \hookrightarrow \textrm{Dom}(A_{\alpha})  \hookrightarrow H^{\frac{\alpha}{2}}(\Omega)$. The proof is complete.
 \end{proof}

\subsection{The local and global solution of problem (\ref{4.1})}

The existence proof of the local solution is a standard contraction argument. With numbers $T>0$ and $R>0$ to be fixed below, in the Banach space $X=C^0([0,\,T],\,H_0^{\frac{\alpha}{2}}(\Omega))$, we consider the closed set
$$
S=\{u\in X:~\|u-u_0\|_X\le R\}.
$$
It follows that the map
\begin{equation}\label{solsemi}
u=\textrm{e}^{tA_{\alpha}}u_0-\int_0^t \textrm{e}^{(t-\tau)A_{\alpha}}(f(u(\tau))-g(x))\intd \tau:=\Theta(u)
\end{equation}
is contraction from $S$ into itself.

Since $\textrm{e}^{tA_{\alpha}}$ is a strongly continuous semigroup, we can choose $T_1$ such that
$\|\textrm{e}^{tA_{\alpha}}u_0-u_0\|_{H^{\frac{\alpha}{2}}(D)}\le R/2$ for $t\in[0,\,T_1]$. Denote $F(x,u)=-f(u)+g(x)$. If $u\in X$,
because $f$ is Lipschitz continuous from bounded subsets of $L^{6}(\Omega)$ to $L^{2}(\Omega)$, then
we have a bound $\|F(x,u)\|_X\le K_3 R$, where $K_3>0$ is a constant. Thus, using Lemmas $\ref{semiest}$ and $\ref{fracembe}$, we have
\begin{eqnarray}\label{intoS}
&&\int_0^t \|\textrm{e}^{(t-\tau)A_{\alpha}}(F(\tau,u(\tau)))\|_{H^{\frac{\alpha}{2}}(\Omega)}\intd \tau \nonumber \\
&\le& C\int_0^t \|\textrm{e}^{(t-\tau)A_{\alpha}}(F(\tau,u(\tau)))\|_{L^2(\Omega)}\intd \tau +C\int_0^t \|A^{\beta}_{\alpha}\textrm{e}^{(t-\tau)A_{\alpha}}(F(\tau,u(\tau)))\|_{L^2(\Omega)}\intd \tau \nonumber \\
&\le& C\int_0^t\|\textrm{e}^{(t-\tau)A_{\alpha}}\|_{L^2(\Omega)}\|F(\tau,u(\tau))\|_{L^2(\Omega)}\intd \tau +C\int_0^t \|A^{\beta}_{\alpha}\textrm{e}^{(t-\tau)A_{\alpha}}\|_{L^2(\Omega)}\|F(\tau,u(\tau))\|_{L^2(\Omega)}\intd \tau \nonumber \\
&\le& CK_3R\int_0^t\textrm{e}^{-\delta(t-\tau)}\intd \tau+CK_3R\int_0^t\frac{\textrm{e}^{-\delta(t-\tau)}}{(t-\tau)^{\beta}}\intd \tau \nonumber \\
&\le& \frac{C K_3R}{\delta}(1-\textrm{e}^{-\delta T_2})+\frac{C K_3R}{1-\beta}T_2^{1-\beta},
\end{eqnarray}
where $\frac12<\beta<1$. If we pick up $T_2\le T_1$ small enough, then such that $\frac{C K_3R}{\delta}(1-\textrm{e}^{-\delta T_2})+\frac{C K_3R}{1-\beta}T_2^{1-\beta}\le R/2$ for $t\in[0,\,T_2]$. Therefore
$\Theta:~S\rightarrow S$ when provided $T\le T_2$.

To arrange that $\Theta$ be a contraction mapping, we also use the Lipschitz continuous properties of $F(\tau,u(\tau))$
for $u,\,\bar{u}\in X$. Hence, for $t\in [0,T_2]$, through Lemmas $\ref{semiest}$ and $\ref{embetoHs}$, we have
\begin{eqnarray}
&&\|\Theta(u)-\Theta(\bar{u})\|_{H^{\frac{\alpha}{2}}(\Omega)}\le\int_0^t \|\textrm{e}^{(t-\tau)A_{\alpha}}(F(\tau,u(\tau))-F(\tau,\bar{u}(\tau)))\|_{H^{\frac{\alpha}{2}}(\Omega)}\intd \tau \nonumber \\
&\le& C\int_0^t \|\textrm{e}^{(t-\tau)A_{\alpha}}(F(\tau,u(\tau))-F(\tau,\bar{u}(\tau)))\|_{L^2(\Omega)}\intd \tau  \nonumber \\
&&+C\int_0^t \|A^{\beta}_{\alpha}\textrm{e}^{(t-\tau)A_{\alpha}}(F(\tau,u(\tau))-F(\tau,\bar{u}(\tau)))\|_{L^2(\Omega)}\intd \tau \nonumber \\
&\le& Cl_f(R)\int_0^t\|\textrm{e}^{(t-\tau)A_{\alpha}}\|_{L^2(\Omega)}\|u-\bar{u}\|_{L^2(\Omega)}\intd \tau \nonumber \\
&&\!+C l_f(R)\!\!\!\int_0^t \|A^{\beta}_{\alpha}\textrm{e}^{(t-\tau)A_{\alpha}}\|_{L^2(\Omega)}\|u-\bar{u}\|_{L^2(\Omega)}\intd \tau \nonumber \\
&\le& \frac{C l_f(R)}{\delta}(1-\textrm{e}^{-\delta T})\|u-\bar{u}\|_{H^{\frac{\alpha}{2}}(\Omega)}
+\frac{C l_f(R)}{1-\beta}T^{1-\beta}\|u-\bar{u}\|_{H^{\frac{\alpha}{2}}(\Omega)},
\end{eqnarray}
where $l_f(R)$ denote the Lipschitz constant and $\frac12<\beta<1$; now if $T\le T_2$ is choosen small enough, then
we get $\|\Theta(u)-\Theta(\bar{u})\|_X\leq L\|u-\bar{u}\|_X$, $L<1$, making $\Phi$ a contraction mapping from $S$ into itself.
Thus $\Phi$ has a unique fixed point $u$ in $S$, solving ($\ref{solsemi}$). We have proved the following result
\begin{theorem}\label{localsol}
If $f$ is Lipschitz continuous locally, then problem (\ref{4.1})
has a unique solution $u\in C^0([0,\,T],\,H_0^{\alpha/2}(\Omega))$, where $T>0$ is chosen above.
\end{theorem}

Following, we prove the global solution to the nonlocal semi-linear equations basing on the result of
the local existence.

\begin{theorem} \label{semithm}
Let $g$ belongs to $ L^2(\Omega)$ and Assumption $III$ hold. Then the solution of problem (\ref{4.1}) exists globally in the space $C^0(\mathbb{R}^+,\,H_0^{\frac{\alpha}{2}}(\Omega))$, $\mathbb{R}^+=\{t\in\mathbb{R}~|~t\geq0\}$.
\end{theorem}
\begin{proof}
It's enough to prove that $\displaystyle{\sup_{0\le t<+\infty}\|u(x,\,t)\|_{H^{\frac{\alpha}{2}}_0(\Omega)}<+\infty}$ by standard energy estimates.

Multiplying $u$ on the both side of the Equation (\ref{4.1}) and integrating over $\Omega$, we have
\begin{eqnarray}
\frac12\frac{\intd}{\intd t}|u|^2+\int_\Omega f(u)u \intd x&=&\int_\Omega uA_{\alpha}u\intd x+\int_\Omega g(x)u \intd x \nonumber \\
&\le& -\int_\Omega\int_\Omega (\mathcal{D}^{\ast}(u))^2\intd y \intd x+C\int_\Omega u^2 \intd x+C|g(x)|^2 .
\end{eqnarray}
Using the nonlocal Poincar\'e inequality $|u|^2\le C\int_\Omega\int_\Omega (\mathcal{D}^{\ast}(u))^2\intd y \intd x$, we get
\begin{equation}\label{ee4.11}
\frac{\intd}{\intd t}|u|^2+C|u|^2+C\int_\Omega|u|^pdx+2\int_\Omega\int_\Omega (\mathcal{D}^{\ast}(u))^2\intd y \intd x \le C|g(x)|^2+C|\Omega|,
\end{equation}where $|\Omega|$ denote the measure of $\Omega$.
Hence we have
$$\frac{\intd}{\intd t}|u|^2+C|u|^2\le C|g(x)|^2+C|\Omega|.$$
By uniform Gronwall inequality, we have
\begin{equation}\label{L2estsemi1}
|u|^2 \le |u_0|^2\textrm{e}^{-Ct}+C|g(x)|^2+C|\Omega|.
\end{equation}
which implies $\displaystyle{\sup_{0\le t<+\infty}|u(x,\,t)|<+\infty}$.

Integrating (\ref{ee4.11}) between $t$ and $t+1$, we get
\begin{align*}
|u(t+1)|^2-|u(t)|^2+C\int_t^{t+1}|u(\tau)|^2 \intd \tau &+\int_t^{t+1}\int_\Omega|u(\tau)|^2d\tau\nonumber\\&+2\int_t^{t+1}\int_\Omega\int_\Omega (\mathcal{D}^{\ast}(u(\tau)))^2\intd y \intd x \intd\tau \le C|g(x)|^2+C|\Omega|.
\end{align*}
By ($\ref{L2estsemi1}$), we obtain
\begin{equation}\label{preHsestsemi1}
\int_t^{t+1}\int_\Omega\int_\Omega (\mathcal{D}^{\ast}(u(\tau)))^2\intd y \intd x \intd\tau \le C|g(x)|^2+\frac12|u(t)|^2+C|\Omega|\le C(|g(x)|^2,|\Omega|,\,|u_0|^2).
\end{equation}

At the same time, let $F(s) = \int^s_0f(s)ds$. By Assumption $III$, we obtain
\begin{equation}
C|s|^{p} - C|\Omega| \leq F(s)\leq C|s|^{p} +C|\Omega|.
\end{equation}
Therefore
\begin{equation}\label{ee4.14}
C\int_{\Omega}|u|^p - C|\Omega| \leq
\int_{\Omega}F(u)dx\leq C\int_{\Omega}|u|^p + C|\Omega|.
\end{equation}

Combining with (\ref{ee4.11}) and (\ref{ee4.14}), we have
\begin{equation}\label{ee4.16}
\frac{\intd}{\intd t}|u|^2+C|u|^2+\int_{\Omega}F(u)dx+\int_\Omega\int_\Omega (\mathcal{D}^{\ast}(u))^2\intd y \intd x \le C|g(x)|^2+C|\Omega|.
\end{equation}

Multiplying $\frac{du}{dt}$ on the both side of the equation (\ref{4.1}) and integrating over the domain $\Omega$, we have
\begin{equation}\label{ee4.15}
\frac12|\frac{du}{dt}|^2+\frac{1}{2}\frac{d}{dt}\int_\Omega\int_\Omega (\mathcal{D}^{\ast}(u))^2\intd y \intd x+\frac{d}{dt}\int_\Omega F(u)dx\leq\frac12|g(x)|^2,
\end{equation}
here we also use H$\ddot{o}$lder inequality and Cauchy inequality.

Combining with (\ref{ee4.16}) and (\ref{ee4.15}), we have

\begin{equation}\label{ee4.19}
\frac{d}{dt}(\int_\Omega\int_\Omega (\mathcal{D}^{\ast}(u))^2\intd y \intd x+\int_\Omega F(u)dx)+C(\int_\Omega\int_\Omega (\mathcal{D}^{\ast}(u))^2\intd y \intd x+\int_\Omega F(u)dx)\leq C|g(x)|^2+C|\Omega|.
\end{equation}

Applying uniform Gronwall inequality, we deduce from (\ref{ee4.19}) that

\begin{equation}
\int_\Omega\int_\Omega (\mathcal{D}^{\ast}(u))^2\intd y \intd x+\int_\Omega F(u)dx\leq e^{-Ct}(\int_\Omega\int_\Omega (\mathcal{D}^{\ast}(u_0))^2\intd y \intd x+\int_\Omega F(u_0)dx)+C|g(x)|^2+C|\Omega|,
\end{equation}which implies $\displaystyle{\sup_{0\le t<+\infty}\|u(x,\,t)\|_{H^{\frac{\alpha}{2}}_0(\Omega)}<+\infty}$.
The proof is complete.
\end{proof}

By Theorem \ref{semithm}, the semigroup $T(t)$ corresponding to problem (\ref{4.1}) can be defined by $T(t):u_0\longrightarrow u(t)$, $u(t)$ is the solution of problem (\ref{4.1}). Also, the semigroup $T(t)$ possesses a bounded absorbing set $B_0$, namely, for any $B\subset H^{\frac{\alpha}{2}}_0(\Omega)$, there exists a time $t_0=t_0(B)$ such that when $t\geq t_0$, $T(t)B\subset B_0$.

\subsection{Proof of Theorem \ref{t4.1}}

Now we ready to prove Theorem \ref{t4.1}.

{\bf Proof of Theorem \ref{t4.1}} Thanks to an embedding theorem, we obtain the well-posedness of problem (\ref{4.1}) in space $H^{\frac{\alpha}{2}}(\Omega)$. By Lemma \ref{l3.1}, the assumptions $I$ and $II$ are all satisfied. Hence by Theorem \ref{t2.1}, next we only need to verify spectral gap condition which similar to the proof of Theorem \ref{t3.2}, we also decompose it into three cases:

{\bf Case 1, $\alpha=1$:} $\lambda_{n+1}-\lambda_{n}=\frac{1}{2}+o(\frac{1}{n}).$ While $n$ large enough, we have $\lambda_{n+1}-\lambda_{n}\geq\frac{1}{2}.$ Thus the Lipschitz constant of nonlinearity do not bigger than $\frac{1}{2}$, the spectral gap condition be satisfied.

Supposing that $G(n)=\lambda_{n+1}-\lambda_{n}= (\frac{n+1}{2}-\frac{2-\alpha}{8})^\alpha-(\frac{n}{2}-\frac{2-\alpha}{8})^\alpha+o(\frac{1}{n})$. Then
\begin{equation}\label{e4.12}
G^{'}(n) = \frac{\alpha}{2}\Big[\big(\frac{n+1}{2}-\frac{2-\alpha}{8}\big)^{\alpha-1}-\big(\frac{n}{2}-\frac{2-\alpha}{8}\big)^{\alpha-1}\Big]+\gamma o(\frac{1}{n^2}).
\end{equation}
Our aim is to obtain $G^{'}(n)>0$, the spectral gap becomes larger and larger, hence for $2l_f$, we can find suitable $N$ such that spectral gap condition hold.

{\bf Case 2, $1<\alpha<2$:} In this case, the term $\frac{\alpha}{2}[(\frac{n+1}{2}-\frac{2-\alpha}{8})^{\alpha-1}-(\frac{n}{2}-\frac{2-\alpha}{8})^{\alpha-1}]>0$, and $\gamma o(\frac{1}{n^2})$ becomes very small when $n$ large enough whenever the sign of $\gamma$. Hence we can choose a big enough $N$, such that $G^{'}(N)>0$, so the spectral gap condition hold.

{\bf Case 3, $0<\alpha<1$:} From (\ref{e4.12}), we know that in this situation $\frac{\alpha}{2}[(\frac{n+1}{2}-\frac{2-\alpha}{8})^{\alpha-1}-(\frac{n}{2}-\frac{2-\alpha}{8})^{\alpha-1}]<0.$ While $n$ large enough, we have $G^{'}(n)<0$. Thus we can not fine suitable $N$ such that the spectral gap condition hold.

In conclusion, we see that in cases 1 and 2, there exists inertial manifold, but in case 3, there is not exists inertial manifold for problem (\ref{4.1}). The proof is complete.\hfill$\Box$

\section{Asymptotic approximation of inertial manifold when normal diffusion   is sufficiently small}
\label{approx}

In this section, we   approximate the inertial manifold  $\mathcal {M}_\varepsilon$,  when the normal diffusion $\varepsilon$ is sufficiently small \cite{11}. We will see the relationship between $\mathcal {M}_\varepsilon$ and $\mathcal {M}_0$, as normal diffusion $\varepsilon$ convergent to $0^+$.

  We know, from section 3, that the inertial manifold $\mathcal {M}_\varepsilon$ is the graph of a Lipschitz continuous mapping
\begin{align}\label{e5.1}
\Phi^\varepsilon(p)=(I-P_N)\varphi(p)(0)=-\int^0_{-\infty}e^{sA_\varepsilon}(I-P_N)(f(\varphi(s))-g(x))ds,~\forall~p\in P_N H,
\end{align}
where $A_\varepsilon=- \varepsilon\Delta  +(-\Delta)^{\frac{\alpha}{2}}.$
That is,
$$\mathcal {M}_\varepsilon=\{p+\Phi^\varepsilon(p)~|~p\in P_N H\}.$$

Similarly, from Section 4, the inertial manifold $\mathcal {M}_0$ is the graph of a Lipschitz continuous mapping
\begin{align}\label{e5.2}
\Phi^0(p)=(I-P_N)\varphi(p)(0)=-\int^0_{-\infty}e^{-sA_\alpha}(I-P_N)(f(\varphi(s))-g(x))ds,~\forall~p\in P_N H,
\end{align}
where $A_\alpha$ is the nonlocal operator. That is,
$$\mathcal {M}_0=\{p+\Phi^0(p)~|~p\in P_N H\}.$$

We expand $\Phi^\varepsilon(p)$ as follows. For $p\in P_N H$, set
\begin{equation}\label{e5.3}
\Phi^\varepsilon(p)=\Phi^0(p)+\varepsilon\Phi^1(p)+\varepsilon^2\Phi^2(p)+\cdot\cdot\cdot+\varepsilon^k\Phi^k(p)+\cdot\cdot\cdot.
\end{equation}
We   write the solution of problem (\ref{3.1}) in the form
\begin{equation}\label{e5.4}
u(t)=u_0(t)+\varepsilon u_1(t)+\varepsilon^2u_2(t)+\cdot\cdot\cdot+\varepsilon^ku_k(t)+\cdot\cdot\cdot,
\end{equation}
with the initial condition
\begin{equation}
u(0)=p+\Phi^\varepsilon(p)=p+\Phi^0(p)+\varepsilon\Phi^1(p)+\varepsilon^2\Phi^2(p)+\cdot\cdot\cdot+\varepsilon^k\Phi^k(p)+\cdot\cdot\cdot.
\end{equation}

At   $\varepsilon=0$, we expand $f(u(t))$ (which depends on $\varepsilon$)  by Taylor expansion,
\begin{equation}\label{e5.6}
f(u)=f(u_0(t))+\varepsilon f^{'}(u_0(t))u_1(t)+\frac{f^{''}(u_0(t))}{2!}u_2(t)\varepsilon^2+\cdot\cdot\cdot+\frac{f^{(k)}(u_0(t))}{k!}u_k(t)\varepsilon^k+\cdot\cdot\cdot,
\end{equation}where $f^{(k)}$ denote the $kth$ Fr\'{e}chet derivative of $f$.

Substituting (\ref{e5.4}) and (\ref{e5.6}) into problem (\ref{3.1}), we have

\begin{align}\label{5.7}
\begin{cases}
\dfrac{du_0(t)}{dt} -A_\alpha u_0(t)+f(u_0(t))=g(x)
,&in\  \Omega.\\
u_0(t)\big|_{\Omega^c}=0.\\
u_0(x,0)=p+\Phi^0(p),
\end{cases}
\end{align}

\begin{align}\label{5.8}
\begin{cases}
\dfrac{du_1(t)}{dt} +A_\varepsilon u_1(t)+f^{'}(u_0(t))u_1(t)=0
,&in\  \Omega.\\
u_1(t)\big|_{\Omega^c}=0.\\
u_1(x,0)=\Phi^{1}(p).
\end{cases}
\end{align}
$$\cdot\cdot\cdot\cdot\cdot$$
\begin{align}
\begin{cases}
\dfrac{du_k(t)}{dt} +A_\varepsilon u_k(t)+\frac{f^{(k)}(u_0(t))}{k!}u_k(t)=0
,&in\  \Omega.\\
u_k(t)\big|_{\Omega^c}=0.\\
u_k(x,0)=\Phi^{k}(p).
\end{cases}
\end{align}
and so on.

Solving the above problems, we obtain
\begin{equation}
u_0(t)=e^{A_\alpha t}u_0(0) - \int^t_0e^{A_\alpha(t-s)}(f(u_0(s))-g(x))ds.
\end{equation}
\begin{equation}
u_1(t)=e^{-A_\varepsilon t}u_1(0) - \int^t_0e^{-A_\varepsilon(t-s)}f^{'}(u_0(s))u_1(s)ds.
\end{equation}
$$\cdot\cdot\cdot\cdot\cdot\cdot\cdot$$
\begin{equation}
u_k(t)=e^{-A_\varepsilon t}u_k(0) - \int^t_0e^{-A_\varepsilon(t-s)}\frac{f^{(k)}(u_0(s))}{k!}u_k(s)ds,
\end{equation}
and so forth.

The right hand side of (\ref{e5.1}) can be represented as
\begin{align}\label{e5.13}
&-\int^0_{-\infty}e^{sA_\varepsilon}(I-P_N)(f(\varphi(s))-g(x))ds
\nonumber\\&=-\int^0_{-\infty}e^{-sA_\alpha}(I-P_N)(f(u_0(s))-g(x))ds \nonumber\\&- \varepsilon\int^0_{-\infty}e^{sA_\varepsilon}(I-P_N)f^{'}(u_0(s))u_1(s)ds\nonumber\\&\cdot\cdot\cdot\cdot\cdot\cdot\nonumber
\\&-\varepsilon^k\int^0_{-\infty}e^{sA_\varepsilon}(I-P_N)\frac{f^{(k)}(u_0(s))}{k!}u_k(s)ds\nonumber\\&-\cdot\cdot\cdot\cdot\cdot\nonumber\\&=I_0+\varepsilon I_1+\cdot\cdot\cdot+\varepsilon^k I_k+\cdot\cdot\cdot,
\end{align}
where $I_0=-\int^0_{-\infty}e^{-sA_\alpha}(I-P_N)(f(u_0(s))-g(x))ds$, $I_k=-\int^0_{-\infty}e^{sA_\varepsilon}(I-P_N)\frac{f^{(k)}(u_0(s))}{k!}u_k(s)ds$, $k\geq1.$

By (\ref{e5.3}) and (\ref{e5.13}), we infer from (\ref{e5.1}) that
\begin{equation}
\Phi^0(p)+\varepsilon\Phi^1(p)+\varepsilon^2\Phi^2(p)+\cdot\cdot\cdot+\varepsilon^k\Phi^k(p)+\cdot\cdot\cdot=I_0+\varepsilon I_1+\cdot\cdot\cdot+\varepsilon^k I_k+\cdot\cdot\cdot.
\end{equation}
Matching the powers of $\varepsilon$, we obtain
\begin{equation}\label{e5.15}
\Phi^0(p)=-\int^0_{-\infty}e^{-sA_\alpha}(I-P_N)(f(u_0(s))-g(x))ds,
\end{equation}
and
\begin{equation}\label{e5.16}
\Phi^1(p)=-\int^0_{-\infty}e^{sA_\varepsilon}(I-P_N)f^{'}(u_0(s))u_1(s)ds,
\end{equation}
$$\cdot\cdot\cdot\cdot\cdot$$
\begin{equation}\label{e5.17}
\Phi^k(p)=-\int^0_{-\infty}e^{sA_\varepsilon}(I-P_N)\frac{f^{(k)}(u_0(s))}{k!}u_k(s)ds,
\end{equation}
and so on.

Thus, we see that if the inertial manifold $\mathcal {M}_\varepsilon$ of   (\ref{3.1}) exists and is a graph of a sufficiently smooth function of $\varepsilon$, then $\Phi^0(p)$ and $\Phi^k(p)$, $p\in P_N H$, $k\geq1$ as obtained above   are well defined.

\begin{theorem}\label{t5.1}
Let $\mathcal {M}_\varepsilon$ be the inertial manifold for the system (\ref{3.1}). Assume that the following conditions hold:
\begin{enumerate}
\item[(i)]Nonlinear function $f$ is sufficiently smooth,
\item[(ii)]There exists  an $N\in\mathbb{N}$  such that
$$\lambda_{N+1}-\lambda_N>2l_f.$$
Then for a sufficiently small $\varepsilon$, the inertial manifold $\mathcal {M}_\varepsilon$ can be represented as
$$\mathcal {M}_\varepsilon =\{p+\Phi^0(p)+\varepsilon\Phi^1(p)+\varepsilon^2\Phi^2(p)+\cdot\cdot\cdot+\varepsilon^k\Phi^k(p)+\cdot\cdot\cdot~|~p\in P_N H\},$$
where $\Phi^0(p)$, $\Phi^k(p)$, $k\geq1$, are in (\ref{e5.15}) and (\ref{e5.17}), respectively.
\end{enumerate}
\end{theorem}

\begin{proof}
By the above analysis, we only need to show that (\ref{e5.15}) and (\ref{e5.17}) are well defined. According to section 4, the existence and uniqueness of $\Phi^0(p)$ obviously. Next, we will verify (\ref{e5.17}) well defined.

Thanks to (\ref{e2.6}), (\ref{e5.4}) and (\ref{e5.6}), for problem (\ref{3.1}) we have
\begin{align}
u(t,p)=J(u,p)(t)=&e^{-tA_\varepsilon}p+\int^0_te^{-(t-s)A_\varepsilon}P_N(f(u(s))-g(x))ds\nonumber\\&-\int^t_{-\infty}e^{-(t-s)A_\varepsilon}(I-P_N)(f(u(s))-g(x))ds.
\end{align}
Equating the terms with the same power of $\varepsilon$, we get
\begin{align}\label{e5.19}
u_k(t,p)=J(u_k,p)(t)=&\int^0_te^{-(t-s)A_\varepsilon}P_N\frac{f^{(k)}(u_0(s))}{k!}u_k(s)ds\nonumber\\&-\int^t_{-\infty}e^{-(t-s)A_\varepsilon}(I-P_N)\frac{f^{(k)}(u_0(s))}{k!}u_k(s)ds.
\end{align}

Under further assumptions, for a large enough $N$, such that $\lambda_{N+1}-\lambda_N>2l_f$, one can choose $\sigma$ such that
\begin{equation}
\lambda_N + 2l_fK_1<\sigma<\lambda_{N+1}-2l_fK_1K_2.
\end{equation}
This $\sigma$ is used to define the Banach space
$$\mathcal {F}_\sigma=\{\varphi\in C((-\infty,0],H)~|~\|\varphi\|_{\mathcal {F}_\sigma}=\sup_{t\leq0}e^{\sigma t}|\varphi(t)|<\infty\}.$$

By section 2,  we know that $u_k(t,p)$ is a well defined mapping from $\mathcal {F}_\sigma\times P_N H\longrightarrow\mathcal {F}_\sigma.$
Indeed, for $|\frac{f^{(k)}(u_0(s))}{k!}|\leq 2l_f$, we have
\begin{align*}
|J(u_k,p)(t)|\leq 2l_f\max\{\int^0_t|e^{-(t-s)A_\varepsilon}P_N||u_k(s)|ds,\int^t_{-\infty}|e^{-(t-s)A_\varepsilon}(I-P_N)||u_k(s)|ds\}.
\end{align*}
By exponential dichotomy properties, we have
\begin{align*}
\|J(u_k,p)(t)\|_{\mathcal {F}_\sigma}&\leq 2l_f\sup_{t\leq0}\max\{\int^0_tK_1e^{(\sigma-\lambda_N)(t-s)}ds,\int^t_{-\infty}  K_2e^{(\sigma-\lambda_{N+1})(t-s)}ds\}\|u_k(s)\|_{\mathcal {F}_\sigma}\nonumber\\&\leq\max\{\frac{2l_fK_1}{\sigma-\lambda_N},\frac{2l_fK_2}{\lambda_{N+1}-\sigma}\}\|u_k(s)\|_{\mathcal {F}_\sigma}<\infty.
\end{align*}

Next we prove that $J(u_k,p)(t)$ is a contraction mapping on $\mathcal {F}_\sigma$.
\begin{align}
\|J(u^1_k,p)(t)-J(u^2_k,p)(t)\|_{\mathcal {F}_\sigma}&=\sup_{t\leq0}e^{\sigma t}\{|\int^0_te^{-(t-s)A_\varepsilon}P_N\frac{f^{(k)}(u_0(s))}{k!}(u_k^1(s)-u_k^2(s))ds+\nonumber\\&\int^t_{-\infty}e^{-(t-s)A_\varepsilon}(I-P_N)\frac{f^{(k)}(u_0(s))}{k!}(u_k^1(s)-u_k^2(s))ds|\}\nonumber\\
&\leq
\sup_{t\leq0}e^{\sigma t}\max\{\int^0_t|e^{-(t-s)A_\varepsilon}P_N\frac{f^{(k)}(u_0(s))}{k!}(u_k^1(s)-u_k^2(s))|ds,\nonumber\\&\int^t_{-\infty}|e^{-(t-s)A_\varepsilon}(I-P_N)\frac{f^{(k)}(u_0(s))}{k!}(u_k^1(s)-u_k^2(s))|ds\}.
\end{align}
By exponential dichotomy condition which was presented in section 2, and $|\frac{f^{(k)}(u_0(s))}{k!}|\leq 2l_f$, we have
\begin{align}
&\|J(u^1_k,p)(t)-J(u^2_k,p)(t)\|_{\mathcal {F}_\sigma}\leq \sup_{t\leq0}e^{\sigma t}\max\{\int^0_t|K_1l_fe^{-\lambda_N(t-s)}|(u_k^1(s)-u_k^2(s))|ds,\nonumber\\&\int^t_{-\infty}K_2l_fe^{-\lambda_{N+1}(t-s)}|(u_k^1(s)-u_k^2(s))|ds\}.\nonumber\\
&\leq
\sup_{t\leq0}\|(u_k^1(s)-u_k^2(s))\|_{\mathcal {F}_\sigma}\max\{K_12l_f\int^0_te^{(\sigma-\lambda_N)(t-s)}ds\nonumber\\&,K_2l_f\int^t_{-\infty}e^{(\sigma-\lambda_{N+1})(t-s)}ds\}\nonumber\\
&\leq\max\{\frac{2l_fK_1}{\sigma-\lambda_N},\frac{2l_fK_2}{\lambda_{N+1}-\sigma}\}\|(u_k^1(s)-u_k^2(s))\|_{\mathcal {F}_\sigma}.
\end{align}
Since by spectral condition, $\frac{2l_fK_1}{\sigma-\lambda_N}<1$ and $\frac{2l_fK_2}{\lambda_{N+1}-\sigma}<1$, hence $J(u_k,p)(t)$ is a contraction mapping on $\mathcal {F}_\sigma$. Using the contraction mapping principle, there exists unique $u_k(t,p)$ satisfy (\ref{e5.19}).

Let $\Phi_k(p)=(I-P_N)u_k(p)(0),~\forall~p\in P_N H.$ We have
$$\Phi^k(p)=-\int^0_{-\infty}e^{sA_\varepsilon}(I-P_N)\frac{f^{(k)}(u_0(s))}{k!}u_k(s)ds,$$
which satisfy (\ref{e5.17}). The proof is complete.
\end{proof}

Following, we consider the asymptotic behavior between $\mathcal {M}_\varepsilon$ and $\mathcal {M}_0$ when $\varepsilon\longrightarrow0^+$. We recall that the inertial manifold for problem (\ref{3.1}) is
$$\mathcal {M}_\varepsilon=\{p+\Phi^0(p)+\varepsilon\Phi^1(p)+\varepsilon^2\Phi^2(p)+\cdot\cdot\cdot+\varepsilon^k\Phi^k(p)+\cdot\cdot\cdot~|~p\in P_N H\},$$
where $\Phi^0(p)$, $\Phi^k(p)$ are as (\ref{e5.15}) and (\ref{e5.17}) respectively. The inertial manifold for problem (\ref{4.1}) is
$$\mathcal {M}_0=\{p+\Phi^0(p)~|~p\in P_N H\},$$
where $\Phi^0(p)$ as in (\ref{e5.15}). In proof of Theorem \ref{t5.1}, we know that $u_k(t,p)\in\mathcal {F}_\sigma$, $p\in P_NH$, hence we have $|\Phi^k(p)|<\infty$, $k\geq1$. Thus, we have following main result in this section

\begin{theorem}
 Let $\mathcal {M}_\varepsilon$ and $\mathcal {M}_0$ be inertial manifolds for problems (\ref{3.1}) and (\ref{4.1}) respectively. Then $\mathcal {M}_\varepsilon$ convergence to $\mathcal {M}_0$ in the norm of $H$ as $\varepsilon\longrightarrow0^+,$ that is $~\forall~0<\varepsilon<<1$
 $$dist_H(\mathcal {M}_\varepsilon,\mathcal {M}_0)\leq o(\varepsilon),~as~\varepsilon\longrightarrow0^+,$$
 where $dist_H(\cdot,\cdot)$ denote Hausdorff semi-distance between two sets in space $H$.
 \end{theorem}

\begin{proof}

 Supposing that in construction of $\mathcal {M}_\varepsilon$ and $\mathcal {M}_0$, we choose same initial data $p\in P_N H$. Then by the expression of $\mathcal {M}_\varepsilon$ and $\mathcal {M}_0$ we get
\begin{align*}
dist_H(\mathcal {M}_\varepsilon,\mathcal {M}_0)\leq\varepsilon|\Phi^1(p)|+\varepsilon^2|\Phi^2(p)|+\cdot\cdot\cdot+\varepsilon^k|\Phi^k(p)|+\cdot\cdot\cdot.
\end{align*}
Since $|\Phi^k(p)|<\infty$, $k\geq1$,  we obtain the result. The proof is complete.
\end{proof}


\begin{thebibliography}{99}


\bibitem{Applebaum} D. Applebaum,
  \textit{L\'evy Processes and Stochastic Calculus}, Second Edition,
 Cambridge University Press, Cambridge, 2009.

\bibitem{31}D. Bl\"omker, \textit{Nonhomogeneous Noise and Q-Wiener Processes on Bounded Domains}, \\Stochastic Analysis and Applications, 23(2), 255-273, 2005.

\bibitem{15}L. Caffarelli, J.-M. Roquejoffre, Y. Sire, \textit{Variational problems with free boundaries for the fractional Laplacian}, J. Eur. Math. Soc. 12 (2010), 1151--1179.


\bibitem{16}L. Caffarelli, L. Silvestre, \textit{An extension problem related to the fractional Laplacian}, Comm. Partial Differential Equations, 32(2007), 1245--1260.


\bibitem{18}L. Caffarelli, A. Vasseur, \textit{Drift diffusion equations with fractional diffusion and the quasi-geostrophic equation}, Annals of Mathematics, 171(2010), 1903--1930.


\bibitem{26}Z. Chen, P. Kim, R.  Song, \textit{Heat kernel estimates for Dirichlet fractional Laplacian}, J. European Math. Soc. 12(2010), 1307--1329.


\bibitem{CMN2012} Z. Chen, M. M. Meerschaert, E. Nane, \textit{Space-time fractional diffusion on bounded domains},
J. Math. Anal. Appl. 393(2012), 479--488.

\bibitem{29}R. Cont, P. Tankov, \textit{Financial modelling with jump processes},   Chapman
\& Hall/CRC, Boca Raton, FL, 2004.


\bibitem{22}A. Debussche, R. Temam, \textit{Convergent families of approximate inertial manifolds}, J. Math. Pures Appl. (9), 75(5)(1994), 489--522.



\bibitem{12}E. Di Nezza, G. Palatucci, E. Valdinoci, \textit{Hitchhiker's guide to the fractional Sobolev spaces,} Bull. Sci. Math., 136(2012), 521--573.



\bibitem{20}J. Duan, K. Lu, B. Schmalfuss, \textit{Invariant manifolds for stochastic
partial differential equations}, \emph{ Ann. Probab.}, 31
(4), 2003, 2109-2135.


\bibitem{13} Q. Du, M. Gunzburger, R. B. Lehoucq, K. Zhou, \textit{A nonlocal vector calculus, nonlocal volume-constrained problems, and nonlocal balance laws,}  Math. Models Methods Appl. Sci. 23(93) (2013). DOI: 10.1142/S0218202512500546.

\bibitem{DGLZ2011} Q. Du, M. Gunzburger, R. B. Lehoucq, and K. Zhou, \textit{Analysis and approximation of
nonlocal diffusion problems with volume constraints}, SIAM Rev., 54 No. 4(2012), 667--696.

\bibitem{17}Q. Du, K. Zhou, \textit{Mathematical analysis for the peridynamic nonlocal continuum theory}, Math. Model. Numer. Anal. 45 (2011), 217--234.


\bibitem{25}T. Gao, J.  Duan, X.  Li, R.  Song, \textit{Mean exit time and escape probaliety for dynamical systems driven by L$\acute{e}$vy noise}, arXiv:1201.6015 (January 2012).


\bibitem{5}M. G. Garroni, J. L. Menaldi, \textit{Second order elliptic integro-differential problems}, Chapaman \&Hall/CRC, 2002.




\bibitem{3}J. K. Hale, \textit{Asymptotic behavior of dissipative systems},
Mathematical Surveys and Monogrphs, vol. 25, AMS, Providence, 1988.

\bibitem{4}D. Henry, \textit{Geometric theory of semilinear parabolic equation}, Lecture Notes in mathematics, 840. Springer-Verlag, BerlinNew York, 1981.




\bibitem{21}M. Kwak, \textit{Finite-dimensional inertial forms for the 2D Navier-stokes equations}, Indiana Univ. Math. J. 41 (1992), no. 4, 927-981.


\bibitem{10}M. Kwasnicki, \textit{Eigenvalues of the fractional Laplacian operator in the interval}, Journal of Functional Analysis, 262(5)(2012), 2379--2402.

\bibitem{MS2012} M. M. Meerschaert and A. Sikorskii, \textit{Stochastic Models for Fractional Calculus}, Walter de Gruyter GmbH \& Co. KG., Berlin/Boston, 2012.


\bibitem{MK2004} R. Metzler and J. Klafter, \textit{The restaurant at the end of the random walk: recent developments
in the description of anomalous transport by fractional dynamics}, Journal of Physics A:
Mathematical and General, 37(31):R161, 2004.

\bibitem{19}X. Mora, J. Sola-Morales, \textit{Existence and nonexistence of finite dimensional globally attracting invariant manifolds in semilinear damped wave equation}, Dynamics of infinite-dimensional systems (Lisbon, 1986), 187--210.


\bibitem{23}J. Novo, E. S. Titi, S. Wynne, \textit{Efficient methods using high accuracy approximate inertial manifold}, Numer. Math., 87(3)(2001), 523--554.


\bibitem{8}A. Pazy, \textit{Semigroup of linear operator amd application to partial differential equation}, Spinger-Verlag, Berlin, 1983.


\bibitem{24}J. C. Robinson, \textit{Computing inertial manifolds}, Discrete Contin. Dyn. Syst., 8(4)(2002), 815--833.

\bibitem{2}J. C. Robinson, \textit{Infinite-Dimensional Dynamical Systems:
An Introduction to Dissipative Parabolic PDEs and the Theorem of
Global Attractors}, Cambridge: Cambridge University Press, 2001.

\bibitem{30}O. Savin, E. Valdinoci, \textit{Elliptic PDEs with fibered nonlinearities}, J. Geom. Anal. 19(2), 420--432.

\bibitem{6}G. R. Sell, Y. You, \textit{Dynamics of evolutionary equations}, Appl. Math. Sci., vol, 143, Springer-Verlag, 2002.


\bibitem{14}L. Silvestre, \textit{Regularity of the obstacle problem for a fractional power of the Laplace operator}, Comm. Pure Appl. Math. 60(1)(2007), 67--112.

\bibitem{28}L. Silvestre, V. Vicol, A. Zlato$\check{s}$, \textit{On the loss of continuity for super-critical drift-diffusion equations}, Arch. Rational Meeh. Anal. 207(2013), 845--877.

\bibitem{11}X. Sun, J.  Duan, X.  Li, \textit{An impact of noise on invariant manifolds in nonlinear dynamical system}, Journal of mathematical physics., 51(2010)042702.


\bibitem{1}R. Temam, \textit{Infinite Dimension Dynamical System in
Mechanics and Physics}, 2nd Edition, Springer, New York, 1997.


\bibitem{9}H. Triebel, \textit{Interpolation theory, function spaces, differential oerators}, North-Holland, Amsterdam-New York, 1978.



\bibitem{7}S. Zelik, \textit{Inertial manifolds and finite-dimensional reduction for dissipative PDEs}, arXiv:1303.4457 [math.AP].




\end{thebibliography}
\end{document}